# Translating scientific Latin texts with artificial intelligence: the works of Euler and contemporaries

Sylvio R. Bistafa
sbistafa@usp.br

## Abstract

The major hindrance in the study of earlier scientific literature is the availability of Latin translations into modern languages. This is particular true for the works of Euler who authored about 850 manuscripts and wrote a thousand letters and received back almost two thousand more. The translation of many of these manuscripts, books and letters have been published in various sources over the last two centuries, but many more have not yet appeared. Fortunately, nowadays, artificial intelligence (AI) translation can be used to circumvent the challenges of translating such substantial number of texts. To validate this tool, benchmark tests have been performed to compare the performance of two popular AI translating algorithms, namely Google Translate and ChatGPT. Additional tests were accomplished in translating an excerpt of a 1739 letter from Johann Bernoulli to Euler, where he announces that he was sending Euler the first part of his manuscript *Hydraulica*. Overall, the comparative results show that ChatGPT performed better that Google Translate not only in the benchmark tests but also in the translation of this letter, highlighting the superiority of ChatGPT as a translation tool, catering not only to general Latin practitioners but also proving beneficial for specialized Latin translators.

Keywords: Latin translation, translating scientific Latin tests, translating Latin with Artificial Intelligence (AI), translating Euler's correspondence

## 1. Introduction

The translation of Euler's works has been a significant endeavor undertaken by mathematicians, historians of mathematics, and scholars over the years. It involves translating the scholarly writings of Leonhard Euler from their original languages (typically Latin or French) into other languages. It allows a wider audience to access and study Euler's mathematical ideas. These translations contribute to preserve Euler's legacy, making his contributions available to mathematicians and researchers worldwide.

Andreas Kleinert in "*Leonhardi Euleri Opera Omnia*: Editing the works and correspondence of Leonhard Euler" [1] gives an overview of the history and present state of the edition of Euler's complete works, since its inception in 1907. This history is characterized by many ups and downs, including the effect of two World Wars and by the necessity to find new publishers for political and financial reasons. Nonetheless, the Swiss Euler Commission has succeeded in the publication of *Leonhardi Euleri Opera omnia* that so far, comprises 81 volumes in four series, which is certainly known to most historians of science and mathematics.

Another source of Euler's translated works is The Euler Archive (http://eulerarchive.maa.org/, https://scholarlycommons.pacific.edu/euler-works/), which is an online dynamic library and database providing summaries and access to digitized versions of original publications, and references to available translations and current research. It contains 207 works with completed translations as of 4 Feb 2021. Knowing that Eneström listed about 850 Euler's publications, there is still much material waiting to be translated. To this end, *Euleriana* (https://scholarlycommons.pacific.edu/euleriana/), an open access eJournal focused on Leonhard Euler and Euler-related scholarship, offers a blind peer-review service for translations to appear in the Euler Archive.

The introductory segment above delves into the intricate task of translating the works of Leonhard Euler, shedding light on the immense effort invested by mathematicians, historians of mathematics, and scholars to make his contributions accessible to a wider audience. It underscores the challenges inherent in translating scientific texts from Latin, emphasizing the meticulous care required in the translation process, especially when dealing with texts from Euler's era.

However, a nuanced perspective suggests that while the focus is ostensibly on Euler's works, the original texts under consideration here predominantly emanate from Euler's contemporaries. This dynamic is not only acceptable but adds a layer of complexity and richness to the scholarly landscape. Acknowledging this, it becomes imperative to expand the scope and articulate the broader scholarly imperative surrounding the translation of works from Latin, particularly those authored by Euler's contemporaries.

The intellectual milieu of the 18th century was characterized by a plethora of groundbreaking mathematical contributions from various scholars, often expressed in Latin—a lingua franca of scientific discourse during that era. This additional layer of context not only refines its focus but also pays due homage to the collaborative nature of mathematical progress during the $18^{th}$ century. It acknowledges the interconnectedness of ideas, fostering a more comprehensive understanding of the translation efforts as a collective endeavor aimed at unraveling the intellectual treasures embedded in the Latin scientific texts of Euler and his esteemed peers.

Translating Latin scientific texts presents difficulties such as specialized technical terminology, language evolution over time, abbreviations and symbols, complex sentence structure, cultural context, lack of contextual information, and the need for translator expertise. Translators must possess subject-specific knowledge, understand the historical context, decipher abbreviations and symbols, and preserve the meaning and logical structure of the original text. Extensive research and consultation may be required to overcome these challenges and ensure faithful translations.

To circumvent these challenges, artificial intelligence (AI) is a rather recent new tool that can be used to ease the translator's tasks. Nonetheless, even with the help of this resource, the translator should have a command of the Latin language, including grammar, vocabulary, syntax, and idiomatic expressions. A basic understanding of Latin grammar rules and the ability to decipher complex sentence structures are also essential. Of course, a command of English is necessary to produce a well-written and coherent translation.

## 2. Latin translation with artificial intelligence AI

AI can be used to translate Latin text, although there are challenges specific to working with ancient languages. These challenges include limited training data, the complex grammar and syntax of Latin, the need for contextual understanding, and domain-specific translation requirements.

While AI can assist in translating Latin, especially for simpler or more common texts, the complexity and nuances of the language can still pose difficulties for automated systems. For tasks requiring high accuracy and precision, human experts with a deep understanding of Latin are often involved to ensure the best quality translations.

Google Translate, Microsoft Translator, or DeepL, are examples of widely used AI translation systems. ChatGPT's has also translation capabilities, although it is important to note that ChatGPT is considered not to be specifically trained or optimized for translation tasks like Google Translate or other dedicated translation tools. ChatGPT is a language model trained on a diverse range of texts from the Internet and is designed to generate human-like responses based on the given context. While it can assist with translation by providing translations as part of a conversation, it lacks the specialized training and fine-tuning that dedicated machine translation systems undergo.

On the other hand, whereas Google Translate is a static translator, ChatGPT is very dynamic and designed for interactive conversations. Although ChatGPT is not specifically optimized for translation like Google Translate, it can generate translations by considering the context of the entire input conversation, allowing for more nuanced translations.

An example of interactive conversation with ChatGPT, where it describes the main characteristics of both translating systems is given below.

SB  What are the main characteristics of translating with Google Translate and ChatGPT?

Google Translate and ChatGPT are two different tools that serve distinct purposes and have different characteristics when it comes to translation. Here are the main characteristics of each:

Google Translate:

1. Automatic translation: Google Translate is an online translation service that provides automatic translations between numerous languages. It uses machine learning algorithms and statistical models to generate translations.
2. Vast language coverage: Google Translate supports a wide range of languages, including commonly spoken languages, regional dialects, and even some ancient languages.
3. User-friendly interface: Google Translate offers an easy-to-use web interface and mobile app, allowing users to input text or speech and receive translations quickly.
4. Real-time translation: Google Translate can provide immediate translations in real-time, making it suitable for on-the-go translation needs.
5. Broad user base and constant improvement: With a vast user base, Google Translate benefits from a large amount of user feedback, which helps improve its translations over time.

> **ChatGPT:**
>
> 1. Language generation and interaction: ChatGPT is an AI language model designed for interactive conversations. While it can handle translation tasks, it is not specifically optimized for translation like Google Translate.
> 2. Contextual understanding: ChatGPT can generate translations by taking into account the context of the entire input conversation, allowing for more nuanced translations.
> 3. Customizability and flexibility: ChatGPT can be fine-tuned and customized for specific applications or domains. Developers can train it on specific translation datasets or use it as part of a larger conversational system.
> 4. Versatility: ChatGPT can be used for a wide range of conversational tasks beyond translation, including answering questions, generating text, providing explanations, and more.
> 5. Language proficiency and accuracy: ChatGPT has been trained on a wide variety of text sources and can generate translations with reasonable accuracy. However, it may not match the precision and fluency of specialized translation systems like Google Translate.

To test the capabilities of AI translation, Google Translate and ChatGPT were applied to translate Latin texts with translations already available elsewhere. The strategy was very simple and objective and consisted of presenting the Latin text to these machine translators, where the resulting translations were then analyzed and compared for correctness against the benchmarks.

After the benchmark testing phase, further tests were performed to translate parts of a letter from Johann Bernoulli to Euler available at The Euler Archive. Unfortunately, the letters from Euler to J. Bernoulli are not yet available in this Archive, and, therefore, it is rather frustrating that we do not know Euler's reaction to the matters raised by J. Bernoulli in his letters.

## 3. Benchmark tests

Benchmark tests were performed to compare the performance of Google Translate and ChatGPT, by using two texts excerpts: one from Daniel Bernoulli, *Hydrodynamica* (1738), and another from Euler's E65 *Methodus inveniendi curvas lineas maximi minimive proprietate gaudentes* … The benchmark translations were taken from the works of Ian Bruce, who is considered to be one of the most prolific translators of 18[th] century historical scientific manuscripts, with an impressive number of translations from Latin into English, including Newton's *Principia*, Euler's *Mechanica*, Daniel Bernoulli *Hydrodynamica*, etc. Ian Bruce maintains a regularly updated commented Internet page at https://www.17centurymaths.com/

*Daniel Bernoulli, Hydrodynamica (1738)*

*SECTIO TERTIA.*

*De velocitatibus fluidorum ex vase utcunque formato per lumen qualecunque effluentium.*

*§.1. Priusquam motum aquarum a gravitate propria ortum definire tentemus, ruminabimur quod in Sectione prima §§.18, 19, 20, 21 & 22 a nobis allatum fuit de principiis ad hoc adhibendis.*

*Recordabimur nempe ascensum potentialem Systematis, cujus singulae partes velocitate qualicunque moventur, significare altitudinem verticalem, ad quam centrum gravitatis illius Systematis pervenit, si singulae particulae motu sursum converso sua velocitate, quantum possunt, ascendere intelligantur, & descensum actualem denotare altitudinem verticalem, per quam centrum gravitatis descendit, postquam singulae particulae in quiete fuerant. Tum etiam memores erimus necessaria ascensum potentialem aequalem esse descensui actuali, quando omnis motus in materia substrata haeret, nihilque de eo in materiam insensibilem aut aliam ad systema non pertinentem transit, & denique motum fluidorum talem proxime esse, ut ubique velocitas reciproce sit proportionalis amplitudini vasis respondenti, qua de re suo loco alia quaedam interjiciemus. Nunc convenit examinare sequentem propositionem.*

HYDRODYNAMICS SECTION THREE.

Tr. by Ian Bruce (2014)

Concerning the velocities of fluids flowing from some kind of vessel through an opening of any kind.

§.1. Before we attempt to define the motion of water arising from its own weight, we will ruminate over what had been brought forwards by us from the principles required to be used for this, in the first section, paragraphs §§.18, 19, 20, 21 & 22.

Certainly we will be remembering the potential ascent [i.e. the power or ability to rise to a certain height] of the system, the individual parts of which are moving with some velocity, to indicate the vertical height, to which the centre of gravity of that system would reach, if the individual parts were understood to be able to rise by changing the velocity of their motion by going upwards, as far as they are able, and the actual descent to denote the vertical height, through which the centre of gravity would fall, after the individual particles were at rest [perhaps to return to the initial situation above, depending on the situation]. Then also we will remember that by necessity the potential ascent is equal to the actual descent, when all the motion stays together in the extended material, and nothing moves across from that into unaffected matter or some other matter not pertaining to the system, and then the motion of the fluids to be approximately such, that everywhere the velocity shall be inversely proportional to the corresponding size [i.e. cross-sectional area] of the vessel, concerning which we will add certain other matters in their place. Now it is agreed to examine the following proposition.

HYDRODYNAMICS, THIRD SECTION.

Tr. by ChatGPT (2023)

**1st Response**

On the velocities of fluids flowing through any shaped vessel and any aperture.

§1. Before we attempt to define the motion of water arising from its own gravity, let us reflect on what was brought forth by us in Section One, §§.18, 19, 20, 21 & 22, regarding the principles to be applied to this matter.

**We shall recall, in particular,** the potential ascent of a system, **where** the individual parts are moved with any velocity, **signifying** the vertical height to which the center of gravity of that **system would rise** if the individual particles were understood to ascend, as far as they can, **with their velocities reversed**. **And we shall note that the actual descent indicates the vertical height** through which the center of gravity descends after the individual particles **have come to rest**. **Furthermore, we shall remember** that the necessary potential ascent is equal to the actual descent when **all motion remains in the underlying matter** and **nothing passes into insensible matter** or **other matter unrelated to the system**. **Finally, we shall find** that the motion of fluids is such that the velocity <u>**is approximately proportional to the corresponding size of the vessel**</u>, and we will expound further on this matter in due course. Now, let us examine the following proposition.

An incorrection regarding the translation of … *reciproce sit proportionalis …* was identified (underlined in bold). The user then asked ChatGPT to generate an alternative translation (Regenerated Response).

**Regenerated Response**

On the velocities of fluids flowing through vessels of any shape and through any aperture.

§1. Before attempting to define the motion of water arising from its own gravity, let us recall what was brought forth by us in the first section, §§.18, 19, 20, 21, and 22, regarding the principles to be applied to this subject.

**We shall remember, indeed**, that potential ascent of a system, **whose** individual parts are moved with any velocity, **signifies** the vertical height to which the center of gravity of that **system rises** if the individual particles are understood to ascend, as far as they can, **by a reverse motion with their own velocity**, **and that actual descent denotes the vertical height** through which the center of gravity descends after the individual particles **have been at rest**. **Then we shall also bear in mind** that necessary potential ascent is equal to the actual descent **when all motion adheres to the material substrate** and **none of it passes into insensible matter** or **into matter not pertaining to the system**, **and, finally,** that the motion of fluids is such that the <u>**velocity is inversely proportional to the corresponding amplitude of the vessel**</u>, which matter we shall discuss further in its proper place. Now it is fitting to examine the following proposition.

Firstly, it is observed that ChatGPT tool Regenerated Response provides a successful correction (underlined in bold) to the incorrection that was pointed out. However, it is interesting to note that the Regenerated Response provided different formulations beyond the error that was corrected from the first response. Perhaps it is a matter of taste which one to choose from, since both translations seem equivalent, however, it is felt that the first translation is, in general, technically better that that provided by the Regenerated Response. Nonetheless, providing an alternative translation enriches the capabilities of ChatGPT, allowing alternatives to choose from.

HYDRODYNAMICS SECTION THREE.

Tr. by Google Translate (2023)

**Of the velocities** of fluids flowing from a vessel of any kind formed **by light** of any kind.

§.1. Before we attempt to define the movement of waters arising from their own gravity, let us consider what we have said in the first section §§18, 19, 20, 21 & 22 about the principles to be used for this purpose.

We will remember that the potential ascent of a system, the individual parts of which are moving at any speed, is to denote the vertical height at which the center of gravity of that system reaches, if each particle is understood to ascend as much as it can by moving upward with its own velocity, and the actual descent is to denote the vertical height through which the center gravity descended, after each particle had been at rest. Then we will also remember that the necessary potential ascent is equal to the actual descent, when all the motion is fixed in the substratum material, and nothing of it passes into insensible matter or other matter not belonging to the system, and finally that **the motion of the fluids is so close** that everywhere the velocity is **reciprocally proportional to the size of the vessel to the respondent**, about which **we will interject some other matters in their place**. It is now appropriate to examine the following statement.

Several problems were identified (in bold) in the translation provided by Google Translate. No further action is provided by this machine translator to correct the translation.

Euler, L. E65 *Methodus inveniendi curvas lineas maximi minimive proprietate gaudentes sive solution problematis isoperimetrici latissimo sensu accepti*.

*Additamentum I.*

*De Curvis Elasticis.*

*2. Sit AB (Fig. 1) lamina Elastica utcunque incurvata; vocetur arcus $AM = s$ et radius osculi curvae $MR = R$; atque, secundum BERNOULLIUM, exprimetur vis potentialis in laminae portione $AM$ contenta hac formula $\int \frac{ds}{RR}$, siquidem lamina sit ubique aequaliter crassa, lata et Elastica atque in statu naturali in directum extensa. Hinc ista erit curvae AM indoles, ut in ea haec expressio omnium minimum obtineat valorem.*

*Quoniam vero in radio osculi $R$ differentialia secundi gradus insunt, ad curvam hac proprietate praeditam determinandam quatuor opus erit conditionibus, id quod cum quaestionis natura apprime convenit. Cum enim Fig.I per datos terminos $A$ et $B$ infinitae laminae Elasticae eaeque eiusdem longitudinis inflecti queant, quaestio non erit determinata, nisi praeter duo puncta $A$ et $B$ simul alia duo puncta seu, quod eodem redit, positio tangentium in punctis extremis $A$ et $B$ praescribatur. Proposita namque lamina Elastica longiori, quam est distantia punctorum $A$ et $B$, ea non solum*

*ita incurvari potest, ut intra terminos A et B contineatur, sed etiam ut eius tangentes in punctis hisce datas teneant directiones.*

*His notatis quaestio de invenienda curvatura laminae Elasticae ex hoc fonte resolvenda ita debet proponi: ut inter omnes curvas eiusdem longitudinis, quae non solum per puncta A et B transeant, sed etiam in his punctis a rectis positione datis tangantur, definiatur ea, in qua sit valor huius expressionis $\int \frac{ds}{RR}$ minimus.*

Tr. by Ian Bruce (2013)

ADDITION IA

ELASTIC CURVES

2. Let AB (Fig. 1) be some curved elastic plate; calling the arc $AM = s$ and the radius of osculation of the curve $MR = R$; and, following Bernoulli, the potential force present in the portion of the plate $AM$ may be expressed by this formula $\int \frac{ds}{RR}$, if indeed the lamina shall be everywhere of equal thickness, width and elasticity and extended sitting in its natural direction. Hence the nature of such a curve $AM$ will be, so that in this, this expression may obtain the smallest value of all.

Because truly differentials of the second order are present in the radius of osculation $R$, there is a need for four conditions for determining a curve with the given property, which agrees with the nature of the first question. For since infinitely many elastic laminas may be able to be curved through the given ends $A$ and $B$ and those of the same length, the question will not be determined unless besides the two points $A$ and $B$ likewise two other points or, what amounts to the same, the position of the tangents at the end points $A$ and $B$ may be prescribed. In so much as with a longer proposed elastic lamina, than is the separation of the points $A$ and $B$, that not only must it be curved thus, so that it may be contained between the ends $A$ and $B$, but also so that its tangents at these given points may maintain given directions.

With these noted the question about finding the curvatures of the elastic laminas being resolved from this source thus must be proposed: so that amongst all the curves of the same length, which not only may pass through the points $A$ and $B$, but also may be tangents at these points with right lines given, that may be defined in which the value of this expression $\int \frac{ds}{RR}$ shall be a minimum.

Tr. by Google Translate (2023)

Appendix I

On Elastic Curves.

2. Let AB (Fig. 1) be an elastic plate curved in any way; it is called the arc $AM = s$ and the radius of the **kiss curve** $MR = R$; and, according to Bernoulli, the potential force contained in the portion $AM$ of the plate will be expressed by this formula $\int \frac{ds}{RR}$, given that the plate is everywhere equally thick, broad and elastic and in its natural state stretched in a straight line. Hence the character of the curve $AM$ will be such that in it this expression obtains the smallest value of all.

But since there are differentials of the second degree in the **radius of the kiss** $R$, to determine the curve endowed with this property four conditions will be needed, which is in perfect agreement with the nature of the question. For since in Fig. 1 **an infinite Elastic plate** of the same length is bent through the given ends $A$ and $B$, the problem will not be determined unless, in addition to the two points $A$ and $B$, there are **two other points at the same time**, or, which comes back to the same thing, the position of the tangents at the extreme points $A$ and $B$ is prescribed. For if the elastic plate is placed longer than the distance between the points $A$ and $B$, it can be bent not only so as to be contained within the boundaries $A$ and $B$, but also so that its tangents at these points hold the given directions.

Having noted these, the problem of finding the curvature of the elastic plate from this source must be resolved in this way: that among all the curves of the same length, which not only pass through the points $A$ and $B$, **but are also touched at these points by the given lines, the one in which it is the smallest value of this expression** $\int \frac{ds}{RR}$.

Several problems were identified (in bold) in the translation provided by Google Translate.

<div style="text-align:center">Tr. by ChatGPT (2023)</div>

<div style="text-align:center">Addendum I.</div>

<div style="text-align:center">On Elastic Curves.</div>

Let there be a curved elastic plate $AB$ (Fig. 1); let the arc $AM$ be denoted as $s$, and the radius of curvature of the curve $MR$ be denoted as $R$. According to Bernoulli, the potential energy in the portion of the plate contained between AM is expressed by this formula $\int \frac{ds}{RR}$, assuming that the plate is uniformly thick, wide, elastic, and in a natural state when extended in a straight line. Hence, the nature of the curve $AM$ will be such that this expression obtains the minimum value.

However, since the second-order differentials are present in the radius of curvature $R$, four conditions are necessary to determine this curve with the given property, which is highly suited to the nature of the question. **Because Fig. I allows an infinitely elastic plate** of the same length to be bent through the given points $A$ and $B$, the question will not be determined unless, in addition to the two points A and B, two other points are simultaneously prescribed, or, in other words, the positions of the tangents at the end points $A$ and $B$ are given. For a given elastic plate longer than the distance between points $A$ and $B$, it can be curved not only to fit within the limits of $A$ and $B$ but also to maintain the given directions of its tangents at these points.

With these notations, the question of finding the curvature of the elastic plate from this source must be posed as follows: among all curves of the same length that not only pass through points $A$ and $B$ but also **touch** the given lines at these points, determine the one in which the value of this expression, $\int \frac{ds}{RR}$, is minimized.

Reproduction extract from ChatGPT where a verbal instruction is given regarding the translation of *tangantur.*

> **SB** Correction: tangantur means tangent
>
> Apologies for the oversight. Let me correct the translation for you.
>
> With these notations, the question of finding the curvature of the elastic plate from this source must be posed as follows: among all curves of the same length that not only pass through points A and B but also are tangent to the given lines at these points, determine the one in which the value of this expression, ds * ∫(ds/RR), is minimized.
>
> **SB** Please, try again: Cum enim Fig.I per datos terminos A\ et B infinitae laminae Elasticae eaeque eiusdem longitudinis inflecti queant
>
> Apologies for the confusion. Here is the corrected translation:
>
> Since, in Figure I, an infinitely elastic plate of the same length can be bent through the given points A and B
>
> **SB** Still not correct:
>
> Apologies for the repeated misunderstanding. Here is the corrected translation:
>
> Given that in Figure I an infinitely elastic plate of the same length can be bent through the given points A and B

ChatGPT repeatedly translates "… *infinitae laminae Elasticae…*" as "… an infinitely elastic plate …", where it should have to be translated as "… an infinity [number] of elastic plates …". Therefore, in this case, a human action is necessary, and a manual correction has to be applied.

To test the capabilities of ChatGPT in translating notation heavy text of mathematical expressions, parts of Section 3 of *De Curvis Elasticis*, extracted from Ian Bruce's modern reproduction of the original Latin text, was copied and pasted into ChatGPT dialog box. The whole interactive conversation follows.

# Translate Latin to English.

December 29, 2023

**Anonymous**
Please, translate the following Latin text into English and include the mathematical expressions in full.

**ChatGPT**
I'm here to help! Please provide the Latin text and any mathematical expressions you'd like translated, and I'll do my best to assist you.

3. Quia solutionem ad coordinatas orthogonales (Fig. 2) accommodari convenit, sumatur recta quaecunque AD pro axe, in qua sit abscissa $AP = x$, applicata $PM = y$; ponatur, uti Methodus tradita iubet, $dy = pdx$, $dp = qdx$; erit elementum curvae $Mm = ds = dx\sqrt{(1+pp)}$. Primum ergo, quia curvae, ex quibus quaesita erui debet, isoperimetrae statuuntur, habebitur ista expressio consideranda $\int dx\sqrt{(1+pp)}$, quae cum generali $\int Z dx$ comparata hunc praebet valorem differentialem $\frac{1}{dx} d \cdot \frac{p}{\sqrt{(1+pp)}}$. Deinde cum sit radius osculi

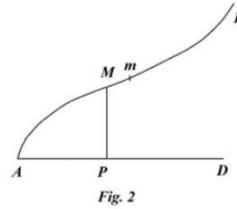

Fig. 2

$$= \frac{dx(1+pp)^{3:2}}{dp} = \frac{(1+pp)^{3:2}}{q} = R,$$

formula $\int \frac{ds}{RR}$, quae minimum esse debet, abit in $\int \frac{qqdx}{(1+pp)^{5:2}}$.

**ChatGPT**
The given Latin text appears to be discussing a solution in orthogonal coordinates, referencing a figure (Fig. 2). The text mentions taking a line AD as an axis with the abscissa AP as $x$. The expression for the curve's element $PM$ is given as $\frac{dy}{dx} = \frac{dp}{q}$, where $Mm$ is the differential $ds$, and $dx$ is compared to $\int Z dx$, resulting in the differential value $\left(1 + \frac{p^2}{1}\right) \frac{1}{p} \frac{d}{dx} \left(1 + \frac{p^2}{1}\right)$.

It then mentions isoperimetric curves, leading to an expression $\int \frac{dx}{\sqrt{1+p^2}}$, which, when compared to $\int Z dx$, yields the differential $\left(1+p^2\right)^{\frac{3}{2}} \frac{d}{dx} \left(1+p^2\right)$. The radius of curvature $R$ is given by $\frac{1}{q}\sqrt{1+p^2}$, and the formula for $ds$ becomes $\int R dx$, which should be minimized and results in $\int \left(1 + \frac{p^2}{1}\right)^{\frac{5}{2}} \frac{q}{p} dx$.

Ian Bruce's (2013) translation follows.

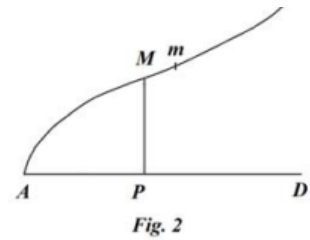

3. Because it will be convenient to adapt the solution to orthogonal coordinates (Fig. 2), some right line $AD$ may be taken for the axis, in which the abscissa shall be $AP = x$, the applied line $PM = y$; there may be put in place, as the method treated orders, $dy = pdx$, $dp = qdx$; an element of the curve will be $Mm = ds = dx\sqrt{(1+pp)}$. Therefore in the first place, because the curves, from which the solution sought must be elicited, are put in place isoperimetric, this expression will have to be considered $\int dx\sqrt{(1+pp)}$, which when compared to the general $\int Zdx$ presents this differential value $\frac{1}{dx} d \cdot \frac{p}{\sqrt{(1+pp)}}$.

Then since the radius of osculation shall be

$$= \frac{dx(1+pp)^{3:2}}{dp} = \frac{(1+pp)^{3:2}}{q} = R,$$

the formula $\int \frac{ds}{RR}$, which must be a minimum, will change into $\int \frac{qqdx}{(1+pp)^{5:2}}$.

It is seen that that ChatGPT omits the ordinate $PM = y$ in the translation, perhaps because ChatGPT is unfamiliar with the technicalities of the Latin noun *applicata* as referred to the vertical orthogonal coordinate (even Ian Bruce's translation did not recognize it as such and translated it as 'applied line'). Regrettably, there are other omissions, and incorrect reproductions of most of the mathematical expressions.

By its turn, Google Translate doesn't even recognize mathematical expressions. In this regard, ChatGPT at least recognizes mathematical formulas and even try to make some sense of them.

It was then asked to ChatGPT to comment on these difficulties. The response follows.

> "The observation you've made highlights some of the strengths and weaknesses of ChatGPT and Google Translate in handling technical or specialized content, particularly mathematical expressions.
>
> ChatGPT's Strengths:
> Recognition of Mathematical Formulas: You correctly pointed out that ChatGPT has the capability to recognize and attempt to make sense of mathematical expressions. This is indeed a valuable feature, especially for users dealing with technical or scientific content.
>
> Contextual Understanding: While there were some omissions and incorrect reproductions in the translation, the fact that ChatGPT attempts to understand and incorporate context, even in technical terms, showcases its versatility. It might

> struggle with specific terms, as you mentioned with "applicata," but the attempt to engage with the content in a meaningful way is noteworthy.
>
> ChatGPT's Limitations:
> Omissions and Incorrect Reproductions: The mentioned omission of the ordinate $PM = y$ and other inaccuracies in mathematical expressions indicate that there are still limitations to the model's understanding, especially in highly technical or specialized domains. It's a reminder that while ChatGPT can handle a broad range of topics, its performance can vary based on the complexity and specificity of the content.
>
> Google Translate's Limitations:
> Lack of Recognition for Mathematical Expressions: As noted, Google Translate doesn't even recognize mathematical expressions, which can be a significant limitation for users dealing with technical or scientific texts. This lack of recognition might result in the loss of crucial information for users seeking accurate translations of such content.
>
> Considerations:
> When dealing with technical content, it's essential to recognize that both ChatGPT and Google Translate may have limitations. It might be beneficial to use these tools as aids rather than sole solutions, especially when precision is critical. Human oversight and expertise remain invaluable, and in cases where technical accuracy is paramount, consulting with a subject matter expert or a professional translator may still be necessary."

While not flawless, ChatGPT outperformed Google Translate in the benchmark tests conducted. The advantage of user interaction and considering the context of the entire conversation enables ChatGPT to deliver more nuanced translations. However, there were cases where a manual correction was required due to an uncertain interpretation.

Here is what ChatGPT says about its own capabilities as compared to Google Translate.

> "Certainly, let's delve into a more detailed commentary on the provided interaction and explore how the capabilities of ChatGPT are described in comparison to Google Translate.
>
> 1. Confidence and Nuance:
> The statement opens by acknowledging that ChatGPT is "not flawless," indicating a level of transparency about its limitations. Despite this, the claim is made that ChatGPT outperformed Google Translate in benchmark tests. This implies a certain confidence in the model's ability to handle language-related tasks better than one of the widely used translation tools.
> The mention of user interaction and context consideration highlights an important aspect of ChatGPT's design. Unlike many traditional language translation tools, ChatGPT is designed for conversation, and its responses can be influenced by the context of the ongoing interaction. This emphasis on context is crucial for providing more nuanced translations, recognizing that language is often complex and context dependent.

2. User Interaction Advantage:
The acknowledgment of the advantage of user interaction is a key point. While Google Translate primarily operates on a one-input, one-output basis, ChatGPT engages in dynamic conversations. This allows it to consider the entire context of a conversation, leading to potentially more accurate and contextually appropriate translations.

3. Manual Correction Instance:
The admission of a specific case where a manual correction was required adds a touch of realism to the statement. It highlights that, despite the overall advantages, there are still instances where the model may struggle or misinterpret certain inputs. This level of honesty is important in managing user expectations.

4. Implications for Future Improvement:
The mention of a specific case requiring manual correction implies that there is ongoing effort to refine and improve the model. Feedback from users on such cases is valuable for the continual development and enhancement of ChatGPT, underlining a commitment to addressing and learning from limitations.

In summary, the statement paints a picture of ChatGPT as a powerful tool with advantages over traditional translation services like Google Translate. It emphasizes the benefits of user interaction and context consideration, while also acknowledging the model's imperfections and the need for occasional manual correction. This combination of confidence, transparency, and a commitment to improvement contributes to a more comprehensive understanding of ChatGPT's capabilities in translation tasks."

## 4. Machine translations of a Latin letter from Johann Bernoulli to Euler

Determining whether one translation is better than another can be subjective and depends on various factors, including context, accuracy, fluency, and the intended audience. Here are some criteria commonly used to characterize the quality of a translation provided by ChatGPT.

Accuracy: A good translation should accurately convey the meaning of the original text without distorting or omitting information.

Fluency: The translation should read smoothly and naturally in the target language, avoiding awkward or stilted expressions.

Naturalness: The translated text should sound natural and idiomatic in the target language, considering cultural nuances and linguistic conventions.

Contextual Appropriateness: Consideration of the context and purpose of the translation is crucial. A good translation aligns with the intended use and audience of the text.

> Consistency: Consistency in terminology and style is important, especially in longer texts or documents.
>
> Cultural Sensitivity: The translator should be mindful of cultural differences to ensure that the translation is culturally appropriate for the target audience.
>
> Faithfulness to the Source: The translation should remain faithful to the tone, style, and intent of the original text.
>
> Grammar and Syntax: Proper grammar and syntax are essential for a high-quality translation. The translated text should adhere to the rules of the target language.
>
> Clarity: The translated text should be clear and easily understandable, avoiding ambiguity or confusion.
>
> Subject Matter Expertise: Depending on the content, a good translation may require knowledge of specialized terminology or subject matter expertise.

It's important to note that evaluating translation quality is often a nuanced process, and what is considered "better" can vary based on the specific requirements of the task and the preferences of the audience. Additionally, feedback from language experts can provide valuable insights into the quality of a translation.

Based on these factors, a comparative quality analysis of the translations provided by ChatGPT and Google Translator was performed on parts of the letter on March 7, 1739, from Johann Bernoulli to Euler available in The Euler Archive. Since different subjects were treated in the same letter by Johann Bernoulli, we have selected only matters related to J. Bernoulli's *Hydraulica*, because of the controversy with his son Danial Bernoulli about his competing publication of *Hydrodymamica*, which Euler seems to have played a mediator over the disputes. Unfortunately, this is a very truncated epistolary exchange because most of the original letters are not still available, particularly the letters from Euler to J. Bernoulli.

J. Bernoulli's work *Hydraulica* is dated 1732 but this is incorrect and was an attempt to obtain priority over his own son Daniel. D. Bernoulli completed his most important work *Hydrodynamica* in 1734 and published it in 1738 at about the same time as Johann published *Hydraulica*.

In the letter on March 7, 1739, J. Bernoulli notifies that he was sending to Euler the first part of his manuscript *Hydraulica*. However, in the preface to the English translations of Hydrodynamics & Hydraulics (Dover Publications, 1968) Hunter Rouse indicates that there is no record that J. Bernoulli had written anything on the subject of hydraulics until some months after D. Bernoulli's treatise was off the press, when he stated in a letter to Euler that he was preparing a manuscript on hydraulics which was already well along (no further details about this letter was given). Nevertheless, J. Bernoulli indicates in the first part that it was written in 1732, a full year ahead of the Russian version of his son's.

**Letter of March 7th, 1739[1]:** J. Bernoulli sends the first part of his hydraulic research. Presentation of his method, which he calls direct, and the advantages it offers over

---

[1] PH Fuss published 14 letters from J. Bernoulli to Euler in his *Correspondance mathématique et physique de quelques célèbres géomètres du XVIIIème siècle*. A copy of this letter is

the one used by Daniel in his Hydrodynamics. Contents of the second part of these same research.

**1st part of the Letter**

> *Exoptatissimae Tuae litterae d. 20. Decembris st. v. mihi traditae sunt atque a me perlectae summa cum voluptate. Ecce! nunc ad Te mitto partem priorem meditationum mea rum hydraulicarum\*), quas tantopere Te desiderare testaris, et vel ideo desideras, quod cognoveris imperfectionem, qua haec doctrina etiamnunc ab aliis tractari soleat, immo, ut candide fateris, Tu ipse frustra omne studium in genuina methodo detegenda collocaveris, invita omni, qua polles, perspicacia. Videbis, originem sequioris successus Scriptorum hydraulicorum ex eo unice venisse, quod nemo hactenus attenderit, partem aliquam finitam virium prementium insumi ad formandum **gurgitem**, quando aqua cogitur ex uno tubo in alium diversae amplitudinis transire, licet **gurges** ipse constare concipiatur ex portiuncula aquae infinite parva. Post pertinacem diutinamque pensitationem animadverti tandem, non sufficere , ut attendatur ad solam illam vim seu pressionem, qua liquor in tubis in motum localem seu progressivum excitetur data cum velocitate, sed praeterea in considerationem trahi debere principium Continuitatis, quo fit ut nulla mutatio in effectibus producendis fiat per saltum, sed successive per gradus infinite parvos, ut in hoc negotio accidit, ubi liquor a velocitate minori ad majorem, vel vicissim a majori ad minorem transire debet; unde omnino necesse est, ut prope transitum, vel ante vel post, concipiatur aliqua portiuncula liquoris , quantumvis parva, cujus stratula infinite parva vel accelerando vel retardando procedant, atque hace portiuncula, inaequabili velocitate gaudens, in **stratulis** est, quam voco **gurgitem**: haec omnia uberius et clarius ex ipso scripto intelliges.*

**ChatGPT translation**

> Your most eagerly anticipated letter of December 20th was delivered to me and I read it with great pleasure. Behold! Now I am sending to you the first part of my hydraulic meditations, which you express such a strong desire for and for which you desire them precisely because you know the imperfection by which this doctrine is still often treated by others. Indeed, as you candidly admit, you yourself have fruitlessly devoted all your efforts to discovering the genuine method, despite your abundant perspicacity. You will see that the origin of the more advanced achievements of hydraulic writers has come uniquely from the fact that nobody until now has observed that a certain finite portion of exerted forces is consumed in forming a **whirlpool** when water is forced to pass from one tube to another of different size, although the **whirlpool** itself is conceived to consist of an infinitely small portion of water. After persistent and prolonged pondering, I finally realized that it is not sufficient to consider only the force or pressure by which the fluid is stimulated to local or progressive motion with a given velocity in the tubes, but also the principle of Continuity must be taken into consideration. This principle dictates that no change in the production of effects occurs abruptly, but rather successively through infinitely small steps, as is the case in this matter where the fluid must pass

---

available in the Euler Archive under the number OO211
http://eulerarchive.maa.org/correspondence/letters/OO0211.pdf

> from a lesser to a greater velocity, or vice versa. Therefore, it is absolutely necessary that, near the transition, either before or after it, a certain small portion of the liquid be conceived, however tiny it may be, whose infinitely **small layers** proceed by accelerating or decelerating, and this portion, enjoying unequal velocity, is what I call a **whirlpool**. You will understand all these things more abundantly and clearly from the very text itself.

**Google Translate translation**

> Your most desired letter d. December 20 v. They were delivered to me and I read them with the greatest pleasure. Behold! I now send you the former part of my meditations on the hydraulic rams*), which you testify that you long for so much, and perhaps because you know the imperfection with which this doctrine is still usually treated by others; you have invested in the method to be discovered, in spite of all the insight you can think of. You will see that the origin of the subsequent success of the hydraulic writings came solely from the fact that no one has hitherto paid attention to the fact that some finite portion of the pressing forces are used to form a **gurgle**, when the water is forced to pass from one tube into another of a different size, although the **gurgle** itself consists of an infinitesimally small portion of water. After persevering and pondering for a long time, I finally noticed that it is not enough to pay attention to that force or pressure alone, by which the liquid in the tubes is excited into a local or progressive motion given the speed, but that the principle of continuity must also be taken into consideration, by which there is no change in the effects produced it may be done by leaps, but successively by infinitely small steps, as happens in this business, where the liquid must pass from a smaller velocity to a greater one, or vice versa from a greater to a smaller one; whence it is absolutely necessary that near the passage, either before or after, some portion of liquid should be conceived, however small, whose infinitesimal **stratums** proceed either by accelerating or decelerating; you will understand more abundantly and clearly from the writing itself.

ChatGPT translates "*gurgitem*" as **whirlpool** (a rapidly rotating mass of water in a river or sea into which objects may be drawn, typically caused by the meeting of conflicting currents), as it refers to the circulatory motion in the zone of separation occupied by the stagnant fluid at an abrupt conduit contraction or expansion[2].

Google Translate translates "*gurgitem*" as **gurgling** (making or characterized by a hollow bubbling sound).

By incorporating the term **whirlpool**, ChatGPT demonstrates its ability to utilize specialized terminology, providing human-like responses based on the given context.

The translation of "*stratulis*" is particularly interesting. It could not be found in a search over the Internet. However, Wiktionary gives a translation for *stratus* as **bed coverings**. By its turn, Google Translate translates it simply as **beds**, whereas ChatGPT gives the proper meaning as **small layers**, also recognizing that ...*lis* at the end of a Latin noun often characterizes a small thing.

---

[2] For more details about this matter see the Preface by Hunter Rouse of the English translations of Hydrodynamics & Hydraulics – Dover Publications, 1968.

> *… Post pertinacem diutinamque pensitationem animadverti tandem, non sufficere, ut attendatur ad solam illam vim seu pressionem, qua liquor in tubis in motum localem seu progressivum excitetur data cum velocitate, sed praeterea in considerationem trahi debere principium Continuitatis, quo fit ut nulla mutatio in effectibus producendis fiat per saltum, sed successive per gradus infinite parvos, …*

ChatGPT translation**:** … After persistent and prolonged pondering, I finally realized that it is not sufficient to consider only the force or pressure by which **the fluid is stimulated to local or progressive motion with a given velocity in the tubes**, but also the principle of Continuity must be taken into consideration. **This principle dictates that no change in the production of effects occurs abruptly, but rather successively through infinitely small steps**, …

Google Translate translation: After persevering and pondering for a long time, I finally noticed that it is not enough to pay attention to that force or pressure alone, by which **the liquid in the tubes is excited into a local or progressive motion given the speed**, but that the principle of continuity must also be taken into consideration, **by which there is no change in the effects produced it may be done by leaps, but successively by infinitely small steps**, …

The highlighted excerpts reveal that the ChatGPT translation is technically more precise than the Google Translate translation which presents some ambiguity and confusion.

**2nd part of the Letter**

> *Videbis etiam methodum hanc directam mirifice conspirare cum indirecta (qua sola usus est Filius meus in sua Hydrodynamica) etenim ambae dant eandem solutionem pro blematum hydraulicorum. Posset autem aliquis mirari, cur, qui ista solvere vult per theoriam virium vivarum, non pariter teneatur rationem habere formandi gurgitis, utpote qui videatur requirere ad sui generationem aliquam partem virium vivarum, aeque ac requiritur partem virium mortuarum; sed causam discriminis explico in scripto meo, monstrans , quantitatem materiae quae componit gurgitem, etsi sit infinite parva, nihilominus opus habere vi finita et determinata pressionis ad acquirendam accelerationem vel re tardationem in stratulis, sive ad id, ut sese gradatim ac commodet ad motum quem liquor jam habet in tubo, in quem ingredi debet. At vero vim vivam quae est in omni materia gurgitis, quippe quae quantitatis est infinite parvae et tantum finitam celeritatem in singulis stratulis habens, oppido patet fore illam vim vivam gurgitis infinite parvam ideoque prorsus incomparabilem cum totali vi viva totins massae aqueae in tubis motae. Hoc ergo notari debuisset a Filio, antequam aggrederetur tractationem Hydraulicae per theoriam conservationis virium vivarum, ne quis scrupulum habere possit, videns negligi considerationem gurgitis, quae in methodo directa citra paralogismum negligi non potest; sed quomodo potuisset hoc praecavere, cum nequidem ideam habuerit naturae gurgitis, quotempore librum suum scripsit. Vides, Vir Clariss., figuras rudi admodum et crassa Minerva esse delineatas, sine ullo ornamento, nedum ad Stereographiae regulas repraesentatas, id sane efficere non potui, si vel maxime voluissem, ob tremorem manuum mearum qui cum aetate continuo ingravescit.*

> *Fortassis dabitur apud Vos aliquis amanuensis qui, Te dirigente, figuras elegantius et majore cum gratia delineare poterit, ita ut ad mentem meam respondeant.*

**ChatGPT translation**

> You will also see that this direct method marvelously aligns with the indirect method (which my son alone used in his Hydrodynamics), as both provide the same solution to hydraulic problems. Someone might wonder why, if one aims to solve these problems through the theory of living forces, they are not equally obliged to consider the formation of whirlpools, as it seems to require a certain portion of living forces for its generation, just as it requires a portion of dead forces. However, I explain the reason for this distinction in my writing, showing that the quantity of matter composing the whirlpool, although infinitely small, nevertheless requires a finite and determinate force of pressure to acquire acceleration or deceleration in its layers, or to gradually adapt to the motion already possessed by the fluid in the tube it is entering. But indeed, the living force present in all the matter of the whirlpool, since it is of infinitely small quantity and has only a finite velocity in each layer, will clearly be infinitely small and, therefore, incomparable to the total living force of the entire mass of water moving in the tubes. Hence, my Son should have taken note of this before embarking on the treatment of Hydraulics through the theory of the conservation of living forces, so that no one could raise any objections upon seeing the neglect of the consideration of whirlpools, which cannot be overlooked in the direct method without falling into fallacy. But how could he have foreseen this when he did not even have the concept of the nature of whirlpools at the time he wrote his book? You see, distinguished Sir, the figures are very rudely and crudely drawn, without any embellishment, let alone adhering to the rules of Stereography. I certainly could not have accomplished this even if I had wished, due to the trembling of my hands, which continually worsens with age. Perhaps you will have a scribe among you who, under your guidance, will be able to draw the figures more elegantly and with greater grace, so that they correspond to my intentions.

The central idea behind J. Bernoulli's criticism on the use of the principle of conservation of living forces by D. Bernoulli's is broadly explained In the Preface of Hydraulics, where J. Bernoulli claimed that the science of hydraulics has not been subjected to the laws of Mechanics, being the developments so far based on experience and theories that were uncertain, and with insufficient foundation, particularly the conservation of living forces, which he recognized as certainly true, and proven by him himself, but that still were not accepted by all philosophers. For him, the true foundations of hydraulics ought to be the principles of mechanics, namely Newton's principles.

**Google Translate translation**

> You will also see that this direct method works wonderfully with the indirect method (which my son used only in his Hydrodynamics), for both give the same solution for hydraulic problems. But one might wonder why he who wishes to solve these things by the theory of living forces, is not equally bound to have a reason for the formation of the gurgling, so that he seems to require for his own generation some part of living forces, just as much as a part of dead forces is required; but I explain the cause of the difference in my writing, pointing out that the amount of

> matter which composes **the gurgling**, even if it is infinitely small, nevertheless needs a finite and determined force of pressure to acquire acceleration or, in fact, retardation **in the beds**, or in order that it gradually and lends itself to the motion which he already has the liquid in the tube into which he must enter. But indeed, the living force which is in all matter of **gurgling**, since it is of infinitesimal quantity and having only a finite speed in each **bed**, it is **clear to the town** that this living force of **gurgling** would be infinitely small and therefore totally incomparable with the total living force of so many masses of water moved in the pipes. This, then, ought to have been noted by **Filius**, before approaching the treatment of hydraulics by the theory of the conservation of living forces, lest anyone should have a scruple, seeing that the consideration of **gurgling** is neglected, which cannot be neglected in the direct method on this side of the paralogism; but how could he have prevented this, since he had no idea of the nature of the **gurgling** at the time he wrote his book. You see, Sir Clariss., that the figures of Minerva are very crudely and crudely drawn, without any ornament, let alone represented according to the rules of stereography, which of course I could not have done, even if I had most wanted to, on account of the trembling of my hands, which grows steadily worse with age. Perhaps some amanuensis will be given to you who, by drawing you, will be able to draw the figures more elegantly and with greater grace, so that they correspond to my mind.

The first thing to be observed is that ChatGPT provides the right expected punctuation to the text, whereas Google Translate provides no punctuation at all.

> *... oppido patet fore illam vim vivam gurgitis infinite parvam ideoque prorsus incomparabilem cum totali vi viva totins massae aqueae in tubis motae ...*

ChatGPT translation: ... will clearly be infinitely small and, therefore, incomparable to the total living force of the entire mass of water moving in the tubes …

Google Translate translation: … **it is clear to the town** that this living force of gurgling would be infinitely small and therefore totally incomparable with the total living force **of so many masses** of water moved in the pipes.

According to the Wiktionary *oppido* it is an Adverb (colloquial in classical texts) meaning very, greatly, much. It is seen that ChatGPT correctly translates it as **infinitely** small; moreover, it gives a more precise mathematical meaning to it (to the infinity) more proper than very, greatly, etc. By its turn, Google Translate improperly translates *oppido* as the dative/ablative singular of *oppidum* meaning town.

> *... prorsus incomparabilem cum totali vi viva totins massae aqueae in tubis motae...*

ChatGPT translation: … and, therefore, incomparable to the total living force of the entire mass of water moving in the tubes ...

Google Translate translation: … and therefore totally incomparable with the total living force of so many masses of water moved in the pipes.

Note the lack of commas before and after "therefore". It is clearly seen that ChatGPT gives a more faithful translation than Google Translate to this excerpt of the letter.

> *... Hoc ergo notari debuisset a Filio, antequam aggrederetur tractationem Hydraulicae per theoriam conservationis virium vivarum, ne quis scrupulum habere possit, videns negligi considerationem gurgitis, quae in methodo directa citra paralogismum negligi non potest; sed quomodo potuisset hoc praecavere, cum nequidem ideam habuerit naturae gurgitis, quotempore librum suum scripsit. ...*

ChatGPT translation: … Hence, my **Son** should have taken note of this before embarking on the treatment of Hydraulics through the theory of the conservation of living forces, **so that no one could raise any objections** upon seeing the neglect of the consideration of **whirlpools**, which cannot be overlooked in the direct method **without falling into fallacy**. But how could he have foreseen this when he did not even have the concept of the nature of whirlpools at the time he wrote his book? …

Google Translate translation: … This, then, ought to have been noted by **Filius**, before approaching the treatment of hydraulics by the theory of the conservation of living forces, **lest anyone should have a scruple**, seeing that the consideration of **gurgling** is neglected, which cannot be neglected in the direct method **on this side of the paralogism**; but how could he have prevented this, since he had no idea of the nature of the **gurgling** at the time he wrote his book. …

*Filio* is translated as **Filius** by Google Translate, whereas it should be translated as **Son** as ChatGPT did.

> *... ne quis scrupulum habere possit* lest anyone …

ChatGPT translation: … so that no one could raise any objections …

Google Translate translation: … lest anyone should have a scruple …

*Scrūpulus*: according to the Wiktionary: m (genitive *scrūpulī*); second declension: (figuratively) anxiety, uneasiness, solicitude, difficulty, **doubt**, **scruple**.

Although the Google Translate translation of *scrupulum* as scruple (a feeling of doubt or hesitation) is not incorrect, the translation given by ChatGPT seems more natural to the targeted language.

### 3rd part of the Letter

> *Ceterum si videro , primam hanc partem hydraulicae meae exercitationis Tibi non displicuisse, transmittam protinus alteram partem , quam interea temporis , dum responsio Tua ad me venerit, absolvam, ut ad mittendum sit parata: Deprehendes, illam adhuc magis esse curiosam, dum ita modifico theoriam meam, ut fere opus non sit idea gurgitis, quem sub alia notione involvo; unde nascitur novum principium hydraulicum, a nemine antea animadversum, cujus auxilio statim pervenio ad motum aquae determinandum fluentis per vasa vel canales, non tantum ex tubis cylindricis conflatos, sed quamcunque figuram, etiam irregularem habentes,*

*aliaque explico phaenomena jucunda et utilia, quae in Physicis quoque suum usum habebunt.*

**ChatGPT translation**

However, if I see that you did not dislike this first part of my hydraulic exercise, I will immediately send you the second part, which I will complete in the meantime while waiting for your response, so that it is ready for transmission. You will discover that it is even more intriguing, as I modify my theory in such a way that **the concept of whirlpools is almost unnecessary**, as I envelop it under a different notion. This gives rise to a new hydraulic principle, unnoticed by anyone before, through which I promptly arrive at determining **the motion of water flowing through vessels or channels, not only formed by cylindrical tubes but also having any shape**, **even irregular ones**. I also explain other pleasant and useful phenomena **that will have their application in physics as well.**

**Google Translate translation**

Moreover, if I see that this first part of my hydraulic exercise is not displeasing to you, I will immediately transmit the second part, which I will complete in the meantime, while your answer comes to me, so that it is ready to be sent. **mine, that the work should not be an idea of gurgling, which I wrap up under another idea**; from which is born a new hydraulic principle, noticed by no one before, with the help of which I immediately arrive at determining **the motion of water flowing through vessels or channels, not only made of cylindrical tubes, but having any shape, even irregular,** and I explain other delightful and useful phenomena, **which in Physics they will also have their use.**

It is seen that the ChatGPT translation is clearer and technically more precise than the Google Translate translation.

## 5. Conclusions

The AI-powered "chat bot" called ChatGPT demonstrated a good proficiency in translating Latin. Despite occasional translation errors, you can verbally instruct ChatGPT to make corrections and provide additional information. Utilizing neural networks, ChatGPT showcases exceptional translating capabilities as a dynamic machine translator, surpassing expectations despite not being explicitly designed or optimized for translation tasks like Google Translate or other dedicated tools. The conducted tests reveal that ChatGPT excels in contextualization, offering alternative vocabulary, syntax, and idiomatic expressions through user interaction. As a result, ChatGPT proves to be a valuable asset for both general Latin translators and those specializing in the language.

However, ChatGPT presents some weakness in handling technical or specialized content, particularly mathematical expressions. The fact that ChatGPT attempts to understand and incorporate context, even in technical terms, showcases its versatility. Google Translate doesn't even recognize mathematical expressions, which can be a significant limitation for users dealing with technical or scientific texts. The difficulties lie not just in language translation but in the interpretation and accurate representation of mathematical concepts, which involve a level of abstraction and precision that may be challenging for machines. Nonetheless, mathematical symbols, equations, and notations often have standardized meanings globally, which often dispense the necessity of translations. This universality is a fundamental characteristic that enables mathematicians worldwide to communicate effectively using a shared symbolic language.

It is, therefore, suggested a cautious approach to applying machine translation in fields where symbolic language plays a critical role, such as mathematics. While technology can be beneficial, the author urges careful consideration and perhaps a reserved stance when it comes to translating mathematical content, especially when it is notation-heavy.

## 6. Closing thoughts on publishing AI-generated translations

In considering the use of AI-generated translations for official publication, it's crucial to strike a balance between the advantages offered by these technologies and the need for human oversight. While AI, such as ChatGPT, can produce impressive results, it's not infallible, and there may be instances where a human touch is necessary for precision and cultural nuances.

I recommend viewing AI-generated translations as a valuable tool in the translation process rather than a standalone solution. Human editors should carefully review and, if necessary, refine the output to ensure accuracy and appropriateness for the target audience.

When presenting AI-generated translations in official publications, transparency is key. It's advisable to explicitly mention that the translation has been generated with the assistance of AI technology. This not only sets clear expectations for the reader but also underscores the role of human judgment in the final product.

In terms of crediting, acknowledgment of the AI tool used (e.g., ChatGPT) and the organization or individual overseeing the human editing process can be included. This dual attribution recognizes the collaborative effort involved in producing the translation.

Ultimately, while AI can significantly expedite the translation process, it should be regarded as a complementary tool in conjunction with human expertise. By combining the strengths of AI and human translators, we can achieve translations that are not only efficient but also uphold the highest standards of accuracy and cultural sensitivity.

## Notes

The translations presented here and some texts provided have been generated with the assistance of AI technology.

OpenAI, the organization behind ChatGPT, recommends providing appropriate context and verifying the information from reliable sources whenever possible, as AI-generated responses should be considered as a tool to assist and augment human knowledge rather than a definitive source of information. Additionally, it's worth noting that the responses generated by ChatGPT are based on its training data up until September 2021, and it does not have access to real-time information or recent events beyond that date.

## Declarations

No funding was received for conducting this study.

The author has no competing interests to declare that are relevant to the content of this article.